\documentclass[12pt]{article}

\usepackage{amssymb}
\usepackage{amsmath,amssymb,amsthm, amscd}
\usepackage{indentfirst}

\theoremstyle{definition}

\title  {   \textbf  {Regular Submanifolds in the Conformal Space ${\mathbb Q}^n_p$  }
 \footnote{  2010 Mathematics Subject Classification: Primary 53A30;
Secondary   53C50.}
  \footnote{ This work is
partially supported by the grant No. D20111007  of Zhongdian NSF
 of Hubei Educational Committee.   }   }

\author {  Changxiong Nie       }

\date{ }

\begin{document}
\maketitle
\noindent\textbf {Abstract.}
 There is a Lorenzian group acting on the conformal space ${\mathbb Q}^n_p$.
 We study the regular submanifolds in the conformal
 space ${\mathbb Q}^n_p$ and construct general submanifold theory
 in the conformal space ${\mathbb Q}^n_p$.
 Finally we give the first
variation formula of the Willmore volume functional of
 submanifolds
 in the conformal space ${\mathbb Q}^n_p$ and classify the conformal
 isotropic submanifolds in the conformal space ${\mathbb Q}^n_p$.

 \medskip

\par\noindent
{\bf {\S} 1. Introduction.}
\par\medskip

   A psudo-riemannian manifold is a manifold with an indefinite metric of index $p(p\geq1)$.
  Such structures arise
  naturally in relativity theory and, more recently,
  string theory. Unlike the considerably more familiar Riemannian manifolds (with the metric index 0),
  Lorentzian manifolds are poorly understood. In this paper
  we study the conformal submanifold geometry in psudo-riemannian space
  forms.

Let $\mathbb R^{N}_{s}$ be the real vector space $\mathbb R^{N}$
with the Lorenzian inner product $\langle,\rangle$ given by
  $$\langle X,Y\rangle
  =\sum_{i=1}^{ N - s } x_iy_i-\sum_{i=N - s + 1 }^{N} x_iy_i,
  \eqno{(1. 1)}$$
  where $X=(x_1, \cdots
x_{ _{N} } ),Y=(y_1, \cdots , y_{ _{N} } )\in\mathbb R^{N}$. We denote by
$C^{n+1}$ the cone in $\mathbb R^{n+2}_{p+1}$ and by $\mathbb Q^n_p$ the
conformal space in $\mathbb R P^{n+1}$:
$$C^{n+1}:=\{X\in\mathbb R^{n+2}_{p+1}|\langle X,X\rangle=0,X\neq0\},\eqno{(1.2)}$$
 $$\mathbb Q^n_p:=\{[X]\in\mathbb R P^{n+1}|\langle X,X\rangle=0\}=
 C^{n+1}/(\mathbb R\backslash\{0\}).\eqno{(1.3)}$$

 Topologically $\mathbb Q^n_p$ is $\mathbb S^{n-p}\times\mathbb S^p/\mathbb
Z_2$, which is endowed by a standard Lorentzian \
metric
 $h=g_{
  _{\mathbb S^{^{n-p}}}    }\oplus(-g_{   _{\mathbb S^{^p}}    })$ and the corresponding
conformal structure $[h]:=\{e^\tau h|\tau\in C^
\infty(\mathbb Q^n_p)\}$.\\

  We define the psudo-riemannian sphere space $\mathbb S^n_p$  and
psudo-riemannian hyperbolic space  $\mathbb H^n_p$ by
$$\mathbb S^n_p=\{u\in \mathbb R^{n+1}_p|(u, u)=1\}, \quad
\mathbb H^n_p=\{u\in \mathbb R^{n+1}_{p+1}|(u, u)=-1\}.   $$
We call $\mathbb R^n_p, \mathbb S^n_p,$ and $\mathbb H^n_p$
psudo-riemannian space forms with index $p$. Denote $\pi=\{[x]\in \mathbb
Q^n|x_{1}=x_{n+2}\},\pi_+=\{[x]\in {\mathbb
Q}^n_1|x_{n+2}=0\},\pi_-=\{[x]\in \mathbb Q^n|x_1=0\}$. Observe the
conformal diffeomorphisms
$$\begin{array}{l}
 \sigma :  \mathbb
R^n_p\rightarrow {\mathbb Q}^n_p\backslash\pi,\quad u\mapsto[(
\frac{( u, u)-1}{2},u, \frac{(u, u)+1}{2})],
   \\
\sigma_+: \mathbb S^n_p\rightarrow  {\mathbb
Q}^n_p\backslash\pi_+,\quad
 u\mapsto [(u,1)],  \\
\sigma_-: \mathbb H^n_p\rightarrow {\mathbb
Q}^n_p\backslash\pi_-,\quad
  u\mapsto [(1,u)]. \\
   \end{array} $$
We may regard    ${\mathbb Q}^n_p$ as the common compactification of
$\mathbb R^n_p$, $\mathbb S^n_p,$ and $ \mathbb H^n_p$, while
 $\mathbb R^n_p$, $\mathbb S^n_p,$ and $\mathbb H^n_p$ as the subsets
 of $\mathbb
 Q^n_p$. Therefore we research the conformal geometry in
 the conformal space $\mathbb
 Q^n_p$ with index p, while it is not necessary to do so in these three
 psudo-riemannian space forms respectively.

  When $p=0$, our analysis in this text can be reduced to the Moebius submanifold geometry in the
  sphere space ({\sl see}  Wang[10]). For more details of Moebius submanifold geometry
  see refs [3, 4, 5, 11, 12], {\sl etc}. Some other results about Lorentz conformal
  geometry see refs.[6-9], {\sl etc}.

This paper is organized as follows. In Section 2 we
  prove the conformal group of the conformal space
 ${\mathbb Q}^n_p$ is $O(n-p+1,p+1)/
 \{\pm I\}$. In Section 3 we
 construct general submanifold theory
 in the conformal space ${\mathbb Q}^n_p$ and give the relationship between conformal invariants and
 isometric
ones for hypersurfaces in Lorentzian space forms. In
Section 4 we give the first
variation formula of the Willmore volume functional of
 regular space-like or time-like submanifolds
 in the conformal space ${\mathbb Q}^n_p$ In
Section 5 we  classify the conformal
 isotropic submanifolds in the conformal space ${\mathbb Q}^n_p$.
 \\

\par\noindent
{\bf {\S} 2. The conformal group of the conformal space
 ${\mathbb Q}^n_p$.}
\par\medskip

 First we introduce

 {\bf Lemma 2.1.}\ Let $\varphi:\mathbf M\rightarrow
 \mathbf M$  be a conformal transformation
 on $m(m>2)$ dimensional psudo-riemannian submanifold $(\mathbf M,g)$,
 {\sl i.e.}, $\varphi$ is a diffeomorphism  and $\varphi^*g=e^{2\tau}g,
 \tau\in C^\infty(\mathbf M)$. If $\mathbf M$ is connected,
 then $\varphi$ is determined by the tangent map $\varphi_*$
 and 1-form $d\tau$ at one fixed point.

 {\sl Proof}\quad  For any point $p\in\mathbf M$, there is a local coordinate $(x^i)$
 around $p$. And
 $(y^i)$ is a local coordinate  around $\varphi(p)$.

 For  psudo-riemannian metric $
 \widetilde{g}=e^{2\tau}g=\varphi^*g$ on $\mathbf M$,
 we denote $\widetilde{\nabla}$ the connection of
 $\widetilde{g}$,
 and $\widetilde{R}$ the curvature tensor,
 $\widetilde{Ric}$ the Ricci curvature tensor.
 Respect to $g$, the corresponding operators are $\nabla,
 R,Ric$, respectively.
 The relation of these operators is as the following equations
 $$
 \widetilde{\nabla}_XY=
 \nabla_XY+X(\tau)Y+Y(\tau)X-g(X,Y)\nabla\tau
 , \ \ \ \ \ \ \ \ \ \ \ \ \ \ \ \ \ \ \eqno{(2.1)}$$
 $$
 \widetilde{R}(X,Y)Z=R(X,Y)Z\ \ \ +\ \ \ g(X,Z)\nabla_Y\nabla\tau
 \ \ \ -\ \ \ g(Y,Z)\nabla_X\nabla\tau\ \ \ \ \ \ \ \ \ \ \ \ \ \ \ \ \ \ \ \ \ \
 \
 $$
 $$
 +\ \ \ [g(X,\nabla\tau)g(Y,Z)-g(Y,\nabla\tau)g(X,Z)]\nabla\tau
  \ \ \ \ \ \ \ \ \ \ \ \ \ \ \ \
 $$
 $$\ \
 +\ \ \ [\nabla_YZ(\tau)+Y(\tau)Z(\tau)-YZ(\tau)
 -g(Y,Z)g(\nabla\tau,\nabla\tau)]X
 $$
 $$\ \ \
 -\ \ \ [\nabla_XZ(\tau)+X(\tau)Z(\tau)-XZ(\tau)
 -g(X,Z)g(\nabla\tau,\nabla\tau)]Y
 ,\eqno{(2.2)}$$
 $$
 \widetilde{R}(X,Y,W,Z)=e^{2\tau}\{
 R(X,Y,W,Z)+g(X,Z)g(W,\nabla_Y\nabla\tau)-g(Y,Z)g(W,\nabla_X\nabla\tau)
 $$
 $$
 +\ \ \ [\ g(X,\nabla\tau)g(Y,Z)-g(Y,\nabla\tau)g(X,Z)\ ]\
 g(W,\nabla\tau) \ \ \ \ \ \ \ \ \ \ \ \ \ \ \ \ \  \ \ \ \ \ \
 $$
 $$
 +\ \ \ [\ \nabla_YZ(\tau)+Y(\tau)Z(\tau)-YZ(\tau)
 -g(Y,Z)g(\nabla\tau,\nabla\tau)\ ]\
 g(W,X)\ \ \ \ \ \ \ \ \ \ \
 $$
 $$
 -\ \ \ [\ \nabla_XZ(\tau)+X(\tau)Z(\tau)-XZ(\tau)
 -g(X,Z)g(\nabla\tau,\nabla\tau)\ ]\
 g(W,Y)
 \},\ \ \ \ \ \eqno{(2.3)}$$
 where $X,Y,Z,W$ are smooth vector fields on $\mathbf M$,
 and $\nabla\tau$ is the gradient of $\tau$ respect to $g$.

 Locally,
 let
 $$
 \nabla_{      \frac{\partial}{\partial x^ { ^ j } }       }
 \frac{\partial}{\partial x^i}  =\sum_k\Gamma^k_{ij}
 \frac{\partial}{\partial x^k}    ,\ \ \ \ \ \ \ \ \ \ \ \
 \nabla_{\frac{\partial}{\partial y^j}}
 \frac{\partial}{\partial y^i}=\sum_k\Gamma^{'k}_{ij}
 \frac{\partial}{\partial y^k},$$
 $$
 g  _  {  _{ij}   }=g(\frac{\partial}{\partial x^i},
 \frac{\partial}{\partial x^ { ^ j } }),\ \
 (g^{ij})=(g  _  {  _{ij}   })^{-1}
 ,\ \
 \varphi_*\frac{\partial}{\partial x^i}=\sum_j
 A^j_i\frac{\partial}{\partial y^j},\ \
 \text d\tau=\sum_iB_idx^i.$$
 First we have
 $$
 g(\varphi_*\frac{\partial}{\partial x^i},
 \varphi_*\frac{\partial}{\partial x^ { ^ j } })\circ\varphi
 =e^{2\tau}g(\frac{\partial}{\partial x^i},
 \frac{\partial}{\partial x^ { ^ j } })
 \eqno{(2.4)}$$
 Acting with $\frac{\partial}{\partial x^k}$ on the both sides of (2.4),
we get
 $$2B_kg(\varphi_*\frac{\partial}{\partial x^i},
 \varphi_*\frac{\partial}{\partial x^ { ^ j } })
 =
 g(\nabla_{
 \varphi_*\frac{\partial}{\partial x^k}
 }
 \varphi_*\frac{\partial}{\partial x^i}
 -\varphi_*
 \nabla_
 {\frac{\partial}{\partial x^k}
  }
 \frac{\partial}{\partial x^i},
 \varphi_*\frac{\partial}{\partial x^ { ^ j } })
 \ \ \ \ \ \ $$
 $$\ \ \ \ \ \ \ \ \ \ \ \ \ \ \ \ \ \ \ \ \ \ \ \ \ \ \ \
     +
 g(\nabla_{
 \varphi_*\frac{\partial}{\partial x^k}
 }
 \varphi_*\frac{\partial}{\partial x^ { ^ j } }
 -\varphi_*
 \nabla_
 {\frac{\partial}{\partial x^k}
  }
 \frac{\partial}{\partial x^ { ^ j } },
 \varphi_*\frac{\partial}{\partial x^i}).$$
 Alternating the positions of $i,j,k$, and by the use of
 $$
 \nabla_{
 \frac{\partial}{\partial x^i}}
 \frac{\partial}{\partial x^ { ^ j } }=
 \nabla_{
 \frac{\partial}{\partial x^ { ^ j } }}
 \frac{\partial}{\partial x^i},\ \
 \nabla_{\varphi_*
 \frac{\partial}{\partial x^i}}
 \varphi_*\frac{\partial}{\partial x^ { ^ j } }=
 \nabla_{\varphi_*
 \frac{\partial}{\partial x^ { ^ j } }}
 \varphi_*\frac{\partial}{\partial x^i},$$
 one will obtain
 $$
 B_ig(
 \varphi_*\frac{\partial}{\partial x^ { ^ j } }
 ,
 \varphi_*\frac{\partial}{\partial x^k})
 +B_jg(
 \varphi_*\frac{\partial}{\partial x^i}
 ,
 \varphi_*\frac{\partial}{\partial x^k})
 -B_kg(
 \varphi_*\frac{\partial}{\partial x^i}
 ,
 \varphi_*\frac{\partial}{\partial x^ { ^ j } })$$
 $$
 =g(\nabla_{
 \varphi_*\frac{\partial}{\partial x^i}
 }
 \varphi_*\frac{\partial}{\partial x^ { ^ j } }
 -\varphi_*
 \nabla_
 {\frac{\partial}{\partial x^i}
  }
 \frac{\partial}{\partial x^ { ^ j } },
 \varphi_*\frac{\partial}{\partial x^k})
 ,\ \ \ \ \ \ \ \ \ \ \ \ \ \ \ \ \ \
 \ \ \ \ \ \ \ \ \ \ \ \ \ \ \ \ \ \ \ \ \
 $$
 and
 $$
 B_kg(\varphi_*\frac{\partial}{\partial x^i},
 \varphi_*\frac{\partial}{\partial x^ { ^ j } })
 =g(\nabla\tau,\frac{\partial}{\partial x^k})
 e^{2\tau}g  _  {  _{ij}   }=g  _  {  _{ij}   }
 g(\varphi_*\nabla\tau,
 \varphi_*\frac{\partial}{\partial x^k}),$$
  where
  $$
  \nabla\tau=\sum_{ij}g^{ij}B_i
  \frac{\partial}{\partial x^ { ^ j } }
  .$$
 Therefore
 $$
 \nabla_{\varphi_*
 \frac{\partial}{\partial x^i}}
 \varphi_*\frac{\partial}{\partial x^ { ^ j } }
 -\varphi_*\nabla_{
 \frac{\partial}{\partial x^i}}
 \frac{\partial}{\partial x^ { ^ j } }=
 B_i\varphi_*\frac{\partial}{\partial x^ { ^ j } }
 +B_j\varphi_*\frac{\partial}{\partial x^k}
 -g  _  {  _{ij}   }\nabla\tau.$$
 We collect the terms of $\frac{\partial}{\partial y^k}$ and get
 $$
 \frac{\partial A^k_j}{\partial x^i}=
 B_iA^k_j+B_jA^k_i+\Gamma^t_{ij}A^k_t
 -g  _  {  _{ij}   }\sum_{st}g^{st}B_sA^k_t-
 \sum_{st}A^s_iA^t_j\Gamma^{'k}_{st}
 .\eqno{(2.5)}$$
 Denote
 $$
 r_{ij}=Ric(\frac{\partial}{\partial x^i},
 \frac{\partial}{\partial x^ { ^ j } }),
 \widetilde{r}_{ij}=\widetilde{Ric}(\frac{\partial}{\partial x^i},
 \frac{\partial}{\partial x^ { ^ j } })
 .$$
 On one hand, form (2.3) we have
 $$
 \widetilde{r}_{ij}=r_{ij}-g  _  {  _{ij}   }\triangle\tau
 +(m-2)[B_iB_j-\frac{\partial B_i}{\partial x^ { ^ j } }
 +\sum_t\Gamma^t_{ij}B_t-g  _  {  _{ij}   }g(\nabla\tau,\nabla\tau)]
 ,\eqno{(2.6)}$$
 where $\triangle$ is the Laplacian respect to $g$.
 On the other hand, we have
 $$
 \widetilde{Ric}(X,Y)=Ric(\varphi_*X,
 \varphi_*Y)\circ\varphi
 ,\eqno{(2.7)}$$
 Therefore
 $$
 \widetilde{r}_{ij}=\widetilde{Ric}(\frac{\partial}{\partial x^i}
 ,\frac{\partial}{\partial x^ { ^ j } })=\sum_{st}A^s_iA^t_jr'_{st}
 ,r'_{st}=Ric(\frac{\partial}{\partial y^s}
 ,\frac{\partial}{\partial y^t})
 .\eqno{(2.8)}$$
 Combining with (2.6) and (2.8), we have
 $$
 \frac{\partial B_j}{\partial x^i}=
 B_iB_j+\sum_t\Gamma^t_{ij}B_t-g  _  {  _{ij}   }
 \sum_{st}g^{st}B_sB_t
 \ \ \ \ \ \ \ \ \ \ \ \ \ \ $$
 $$\ \ \ \ \ \ \ \ \ \ \ \ \ \ \ \ \ \ \ \ \ \ \ \
 +\frac{1}{m-2}(r_{ij}
 -g  _  {  _{ij}   }\triangle\tau-\sum_{st}
 A_i^sA^t_jr'_{st})
 .\eqno{(2.9)}$$
 To the first order ODE (2.5),
 (2.8), one may notice
 $
 A^k_j=\frac{\partial\varphi_k}
 {\partial x^ { ^ j } },B_j=\frac{\partial\tau}{\partial x^ { ^ j } }$.
 If $\mathbf M$ is connected,
 then $\varphi$ is determined by the tangent map $\varphi_*$
 and 1-form $d\tau$ at one fixed point.
$\Box$

 {\bf Theorem 2.1.} Suppose that $\varphi$ is a conformal
 transformation on $\mathbb Q^n_p$,
 $\varphi^*h=e^{2\tau}h$,
 and $x_ { _0 }$ is a fixed point of $\varphi$, then there is $A\in O(n-p+1,p+1)$,
 such that $\varphi=\Phi_A$ and $\Phi_A([X])=[XA]$.

 {\sl Proof}\quad Let  $(\mathbf U,x^i)$ be a coordinate chart around $x_ { _0 }$.
 At point $x_ { _0 }$, denote
    $$
    \frac{\partial\varphi_i}{\partial x^ { ^ j } }\mid_{x_ { _0 }}
    =A^i_j,\ \
    \frac{\partial\tau}{\partial x^ { ^ j } }|_{x_ { _0 }}
    =B_j,\ \
    h_{ij}=h(\frac{\partial }{\partial x^i},\ \
    \frac{\partial }{\partial x^ { ^ j } })|_{x_ { _0 }},\ \
    (h^{ij})=(h_{ij})^{-1}.$$
 Suppose that
 $$
 x_ { _0 }=[u_ { _0 }],\ \ u_ { _0 }=(u_p,u_2)\in\mathbb S^{n-1}\times
 \mathbb S^1\subset\mathbb R^{n+2}_2,\ \
 Ju_ { _0 }=(u_p,-u_2)
 .$$
 And if
 $$
 e_i\in E^\bot_{u_ { _0 }},\ \ \pi_*e_i=
 \frac{\partial }{\partial x^i}|_{x_ { _0 }}
 .$$
 $\{u_ { _0 },Ju_ { _0 },e_p,\cdots,e_n\}$ construct a basis of
 $\mathbb{R}^{n+2}_2$, then there is a orthonormal decomposition of $\mathbb{R}^{n+2}_2$:
 $$
 \mathbb{R}^{n+2}_2=\text{span}\{u_ { _0 },Ju_ { _0 }\}\oplus
\text{span}\{e_p,\cdots,e_n\}
 .$$
 Define linear transformation $A:\mathbb R^{n+2}_2
 \rightarrow\mathbb R^{n+2}_2$ on this basis:
   $$
   A(u_ { _0 })=e^{-\tau(x_ { _0 })}u_ { _0 },\ \
   A(e_i)=e^{-\tau(x_ { _0 })}(\sum_jA^j_ie_j-B_iu_ { _0 })
   ,\eqno{(2.10)}$$
   $$A(Ju_ { _0 })=e^{\tau(x_ { _0 })}Ju_ { _0 }+2e^{-\tau(x_ { _0 })}
   (\sum_{ijk}h^{jk}B_jA_k^ie_i-\sum_{ij}h^{ij}B_iB_ju_ { _0 })
   .\eqno{(2.11)}$$

   First, it is easy to know that$A\in O(n-p+1,p+1)$. In fact,
   it is guaranteed by $
   \sum_{st}A^s_iA^t_jh_{st}=h_{ij}e^{2\tau(x_ { _0 })}$
   (check it on the basis).

   Otherwise, we have
 $$
 \Phi_A(x_ { _0 })=\varphi(x_ { _0 })=x_ { _0 }
 ,\eqno{(2.12)}$$
 $$
 \Phi_{A*}|_{x_ { _0 }}(\frac{\partial}{\partial x^i})
 =\pi_*|_{\mathrm{T}_{x_ { _0 }}\mathbb Q^n_p}
 \circ A\circ ( \pi_*|_{\mathrm{T}_{x_ { _0 }}\mathbb Q^n_p} )^{-1}
 (\frac{\partial}{\partial x^i})
 =
 \sum_jA^j_ie_j=\varphi_*|_{x_ { _0 }}(\frac{\partial}{\partial x^i})
 .\eqno{(2.13)}
 $$

 Suppose that $[u]\in\mathbb Q^n_p$, for any $X,Y\in\mathrm{T}
 _{[u]}\mathbb Q^n_p$, there are $\alpha,\beta\in
 E^\bot_u\subset \mathrm{T}_uC^{n+1}$ such that
 $$
 \pi_*\alpha=X,\ \ \pi_*\beta=Y.$$
 Therefore, from (2.1) we have
 $$
 (\Phi^*_Ah)_{[u]}(X,Y)=(\Phi^*_Ah)_{[u]}(
 \pi_*\alpha,\pi_*\beta)=(\pi^*\circ\Phi^*_Ah)_u(\alpha,\beta)
 \ \ \ \ \ \ \ \ \ \ \ \ \ \ \ \ \ \ \ $$
  $$
  =(A\circ\pi^* h)_u(\alpha,\beta)
  =(\pi^*h)_{A(u)}(\alpha A,\beta A)=\frac{2}{|A(u)|^2}
 \langle \alpha A,\beta A\rangle\ \ \ \ \ \ \ \ \ \ \ \ \ \ \ \ \ \ \ \
 \ \ \
 $$
 $$
 =\frac{|u|^2}{|A(u)|^2}\cdot \frac{2}{|u|^2}\langle \alpha,\beta\rangle
 =\frac{|u|^2}{|A(u)|^2}h_{[u]}(X,Y).
 \ \ \ \ \ \ \ \ \ \ \ \ \ \ \ \ \ \ \ \ \ \ \ \ \  \ \ \ \ \ \ \
  \ \ \ \ \ \  \eqno{(2.14)}$$
 Therefore $\Phi^*_Ah=\frac{|u|^2}{|uA|^2}h$.
 Next we prove that
 $$
 \frac{\partial}{\partial x^i}|_{   _{x_ { _0 }}   }
 (\frac{|u|^2}{|uA|^2})=e^{2\tau(x_ { _0 })}B_i
 .\eqno{(2.15)}$$
 Suppose that there is a local lift of
 $\mathbb Q^n_p$ around $x_ { _0 }\in\mathbb Q^n_p$ such that
 $u:\mathbf U\subset\mathbb Q^n_p\rightarrow
 C^{n+1}$. Then $\pi\circ u=$id, and
 $$
 \frac{\partial u}{\partial x^i}|_{x_ { _0 }}=
 u_*( \frac{\partial }{\partial x^i}|_{x_ { _0 }})=
 u_*\circ\pi_*(e_i)=
 (\pi\circ u)_*(e_i)=e_i
 .\eqno{(2.16)}$$
 Suppose that
 $$
 u=au_ { _0 }+bJu_ { _0 }+\sum_ic^ie_i
 ,$$
 where $a,b,c^i$ are local smooth functions.
 Without difference, we   let
 $$
 a(x_ { _0 })=1,\ b(x_ { _0 })=0,\ c^i(x_ { _0 })=0
 .\eqno{(2.17)}$$
 Using (2.10) and (2.11), we denote $A(u)$ by
 $$
 A(u)=(a-2h(\nabla\tau(x_ { _0 }),\nabla\tau(x_ { _0 }))b-
 \sum_iB_ic^i)u_ { _0 }
 +    be^{\tau(x_ { _0 })} Ju_ { _0 }
 \ \ \ \ \ \ \ \ \ \ \ \ \ \ \ \ \
 $$
 $$ +
 e^{-\tau(x_ { _0 })}  \sum_{ik}  (  2b\sum_jB_jh^{jk}
 +c^k)A^i_k e_i
   :=a'u_ { _0 }+b'Ju_ { _0 }+\sum_ic'^ie_i.
  \eqno{(2.18)}$$
It is easy to check that
 $$
 \frac{\partial a}{\partial x^i}|_{x_ { _0 }}=0,\
 \frac{\partial b}{\partial x^i}|_{x_ { _0 }}=0,\
 \frac{\partial c^j}{\partial x^i}|_{x_ { _0 }}=\delta^j_i
 . \eqno{(2.19)}$$
 Consequently,
 $$
 \frac{\partial  }{\partial x^i}|_{x_ { _0 }}
 ( |u|^2  )= \frac{\partial  }{\partial x^i}|_{x_ { _0 }}
 (2a^2+2b^2+\sum_{jk}c^jc^k\langle e_j,e_k\rangle
 =0. \eqno{(2.20)}$$
 $$
 \frac{\partial  }{\partial x^i}|_{x_ { _0 }}
 ( |A(u)|^2  )= 4\langle \frac{\partial a' }{\partial x^i}|_{x_ { _0 }},
 a'(x_ { _0 })\rangle =-2e^{-2\tau(x_ { _0 })}B_i
 . \eqno{(2.21)}$$
 Thereby
 $$
 \frac{\partial  }{\partial x^i}|_{x_ { _0 }}
 ( \frac{ |u|^2}{|A(u)|^2} )=
 -\frac{|u_ { _0 }|^2   \frac{\partial  }{\partial x^i}|_{x_ { _0 }}
( |A(u)|^2  ) }  { |A(u_ { _0 })|^4   }=
 e^{2\tau(x_ { _0 })}B_i
 . \eqno{(2.22)}$$
  From Lemma 2.1, we have $\Phi_A=\varphi$.$\Box$

  Suppose that  for some fixed point $x_ { _0 }=[(a,b)]\in\mathbb Q^n_p$,
  a conformal transformation
  $\varphi:\mathbb Q^n_p\rightarrow\mathbb Q^n_p$ have
  $$
  \varphi([a,b)]=[(c,d)]
  ,$$
  where
  $$
  (a,b),(c,d)\in\mathbb S^{n-p}\times
 \mathbb S^p
 .$$
 We can certainly find $C\in O(n-p+1 ),D\in O(p+1)$,
 such that $a=cC,b=dD$. That is, $A_p=\text{diag}(C,D)\in
 O(n-p+1,p+1)$ such that $\Phi_{A_p}[(c,d)]=[(a,b)]$.
 Clearly, the conformal transformation
 $\Phi_{A_p}\circ\varphi$ of $\mathbb Q^n_p$ has fixed point $x_ { _0 }$.
 From the above theorem, there is $A\in O(n-p+1,p+1)$ such that $\Phi_{A_p}\circ\varphi
 =\Phi_A$. Thus $\varphi=\Phi_{AA_p^{-1}}$.
 At last, since
 $$
 \Phi:O(n-p+1,p+1)\rightarrow\text{the conformal group of\ }\mathbb Q^n_p
 ,A\mapsto\Phi_A
 $$
 is a epimorphism and ker$(\Phi)=\{\pm I\}$,    we have

 {\bf Theorem 2.2. } \ The conformal group of the conformal space
 ${\mathbb Q}^n_p$ is $O(n-p+1,p+1)/
 \{\pm I\}$.\\

\par\noindent
{\bf {\S} 3. Fundamental equations of submanifolds.}
\par\medskip

 Suppose that $x:\mathbf M\rightarrow{\mathbb Q}^n_p$ is an $m$-dimensional
 psudo-riemannian submanifold. That
is, $x_*(\mathrm{T}\mathbf M)$ is non-degenerated subbundle of
$(\mathrm{T}{\mathbb Q}^n_p,h)$. Let $y:U\rightarrow C^{n+1}$ be a lift
of $x:\mathbf M\rightarrow{\mathbb Q}^n_p$ defined in an open subset
$U$ of $\mathbf M$. We denote by $\Delta$ and $\kappa$ the Laplacian
operator and the normalized scalar curvature of the local non-degenerated
metric $\langle\text dy, \text dy\rangle$. Then we have

{\bf Theorem 3.1.}  On $\mathbf M$ the 2-form $g:= \pm (\langle\Delta y,
\Delta y\rangle-m^2\kappa)\langle\text dy, \text dy\rangle $ is a
  globally defined invariant of $x:\mathbf M\rightarrow {\mathbb Q}^n_p$
  under the Lorentzian group transformations of ${\mathbb Q}^n_p$.

{\sl Proof}\quad  First we can check it out that the expression of
$g$ is invariant to different local lifts. Suppose that
$y:U\rightarrow C^{n+1},\tilde y: \tilde U\rightarrow C^{n+1}$ are
different lifts of $x:\mathbf M\rightarrow{\mathbb Q}^n_p$ defined
in  open subsets $U$ and $\tilde U$ of $\mathbf M$. For the local
non-degenerated metrics $\langle,\rangle_y=\langle\text d  y,\text
d y\rangle$, we denote by $\Delta$ the Laplacian, by $\nabla f$ the
gradient of a function $f$, and by $\kappa$ the normalized scalar
curvatures. And for $\langle\text d\tilde y,\text d\tilde y\rangle$,
we denote by $\tilde\Delta$  the Laplacian, and by $\tilde\kappa$
the normalized scalar curvatures. On $U\cap\tilde U$, we find that
$\tilde y=e^\tau y$, where $\tau$ is local smooth function on
$U\cap\tilde U$. Therefore $\langle\text d\tilde y,\text d\tilde
y\rangle=e^{2\tau}\langle\text d  y,\text d y\rangle,$ and they are
conformal on $U\cap\tilde U$. We have
$$ \widetilde{\omega^j_i}=\omega^j_i+\tau_i\omega^j-\tau^j\omega_i
+\delta_i^j\mathrm{d}\tau,\eqno(3.1)$$
$$e^{2\tau}\widetilde{\Delta}f=\Delta f+(m-2)\langle\nabla\tau,\nabla
f\rangle_y.\eqno(3.2)$$
$$e^{2\tau}\widetilde{\kappa}=\kappa-\frac{2}{m}\Delta\tau
-\frac{ m-2}{m}\langle\nabla\tau,\nabla\tau\rangle_y.\eqno(3.3)$$

It follows that
$$
(\langle\Delta y, \Delta y\rangle-m^2\kappa)\langle\text d  y,\text
d y\rangle
 =
 (\langle\tilde \Delta\tilde  y,\tilde  \Delta\tilde  y\rangle-m^2\tilde \kappa)
 \langle\text d \tilde  y,\text d \tilde y\rangle.\eqno(3.4)$$

 If there is a Lorenzian rotation $T$ acting on $\mathbb Q^n_p$
 and $y:U\rightarrow C^{n+1}$ is a lift of $x:\mathbf M\rightarrow{\mathbb Q}^n_p$ defined
in  open subsets $U$, then
 the submanifold $\tilde x=x\circ T$ must have a local lift like $\tilde y=e^\tau yT$.
Since $T$ perserves the Lorentzian inner product and the dilatation
of the local lift $y$ will not impact the term $(\langle\Delta y,
\Delta y\rangle-m^2\kappa)\langle\text d  y,\text d y\rangle$, the
2-form $g$ is conformally invariant.
$\Box$

 {\bf Definition 3.1.} We call an
$m$-dimensional submanifold $x:\mathbf M\rightarrow {\mathbb Q}^n_p$
a regular submanifold if the 2-form $g:= \pm(\langle\Delta y, \Delta
 y\rangle-m^2\kappa)\langle\text dy, \text dy\rangle$ is non-degenerated.
  And $g$ is called the conformal metric of the regular
submanifold
 $x:\mathbf M\rightarrow {\mathbb Q}^n_p$.

 In this paper we assume that $x:\mathbf M\rightarrow {\mathbb Q}^n_p$
 is a regular  submanifold. Since the metric $g$ is
non-degenerated (we call it the conformal metric), there exists a
unique lift $Y:\mathbf M\rightarrow C^{n+1}$ such that
$g=\langle\text dY,\text dY\rangle$ up to a signature. We call $Y$
the canonical lift of $x$. By taking $y:=Y$ in (3.1) we get
$$\langle\Delta Y, \Delta Y\rangle=m^2\kappa\pm1.\eqno{(3.5)} $$

 Theorem 3.1 implies that

{\bf Theorem 3.2.} Two submanifolds $x,\tilde{x}:\mathbf M\rightarrow
{\mathbb Q}^n_p$
 are
conformal equivalent if and only if there exists $T\in O(n-p+1,p+1)$ such
that $\tilde Y=YT$, where $Y,\tilde Y$ are canonical lifts of
$x,\tilde x$
, respectively .

Let $\{e_1, \cdots , e_m\}$ be a local   basis of $\mathbf
M$ with dual basis $\{\omega^1, \cdots , $ $\omega^m\}$.  Denote
$Y_i=e_i(Y)$. We define
 $$N:=-\frac{1}{m}\Delta Y-\frac{1}{2m^2}\langle\Delta Y,
 \Delta Y\rangle Y, \eqno{(3.6)}$$
then we have
 $$\langle N, Y\rangle=1, \langle N, N\rangle=0,
 \langle N, Y_k\rangle=0,\quad1\leq k\leq m.\eqno{(3.7)}$$
 And we may decompose $\mathbb R^{n+2}_{p+1}$ such that
$$\mathbb R^{n+2}_{p+1}=\text{span}\{Y, N\}\oplus \text{span}
\{Y_1, \cdots , Y_m\}\oplus\mathbb V \eqno{(3.8)}$$ where $\mathbb
V\bot\text{span}\{Y, N, Y_1, \cdots , Y_m\}$.  We call $\mathbb V$
 the conformal normal bundle for $x:\mathbf M\rightarrow {\mathbb Q}^n_p$.
 Let $\{\xi_{m+1}, \cdots , \xi_n\}$ be a local basis for the
bundle $\mathbb V$ over $\mathbf M$. Then $\{Y, N, Y_1, \cdots ,
 Y_m, \xi_{m+1}, \cdots , \xi_n\}$ forms a moving frame in
 $\mathbb R^{n+2}_{p+1}$ along $\mathbf M$. We adopt the conventions on the ranges
of indices in this paper:
  $$1\leq i, j, k, l,r,q\leq m;\quad m+1\leq\alpha, \beta,\gamma,\nu\leq n
  .\eqno{(3.9)}$$
 We may write the
 structure equations as follows
 $$\mathrm{d}Y=\sum_i\omega^iY_i;\quad
 \mathrm{d}N=\sum_i\psi^iY_i+\sum_\alpha\phi^\alpha \xi_\alpha ; \eqno{(3.10)}$$
 $$\mathrm{d}Y_i=-\psi_iY-
 \omega_iN+\sum_j\omega_{i}^{j}Y_j+\sum_\alpha\omega^\alpha_i
 \xi_\alpha ;\eqno{(3.11)}$$
 $$\mathrm{d}\xi_\alpha =
 -\phi_\alpha Y     +  \sum_i\omega_{\alpha}^iY_i+\sum_\beta\omega^\beta_\alpha
 \xi_\alpha , \eqno{(3.12)}$$
 where the coefficients of $\{Y,N,Y_i,\xi_\alpha \}
$ are 1-forms on $\mathbf M$.  It is clear that
 $\mathbb A:=\sum_i\psi_i\otimes\omega^i,
\mathbb
B:=\sum_{i,\alpha}\omega^\alpha_i\otimes\omega^ie_\alpha,\Phi:=\sum_\alpha\phi^\alpha
\xi_\alpha $ are globally defined conformal invariants. If we  denote
$$
\psi_i=\sum_jA_{ij}\omega^j,
\quad\omega^\alpha_i  =   \sum_jB^\alpha_{ij}\omega^j,
\quad  \phi^\alpha  =   \sum_{i  }C^\alpha_i\omega^i,
\eqno{(3.13)}$$ then we can define the covariant derivatives of
these tensors and curvature tensor with respect to conformal metric
$g$:
$$\sum_jC^\alpha_{i, j}\omega^j=dC^\alpha_i-\sum_jC^\alpha_j\omega_{i}^{j}+\sum_\beta
C^\beta_i\omega^\alpha_\beta; \eqno{(3.14)}$$
  $$\sum_kA_{ij, k}\omega^k =  dA_{ij} - \sum_kA_{ik}\omega^k_j
  -  \sum_kA_{kj}\omega^k_i;\eqno{(3.15)}$$
  $$\sum_kB^\alpha_{ij, k}\omega^k=dB^\alpha_{ij}  -  \sum_kB^\alpha_{ik}\omega^k_j
  -  \sum_kB^\alpha_{kj}\omega^k_i+\sum_\beta
B^\beta_{ij}\omega^\alpha_\beta; \eqno{( 3.16)}$$
  $$d\omega^{i}_{j}   +  \sum_k\omega^{i}_{k}\wedge\omega^k_j
  =  \omega^{i}\wedge\psi_{j}  +  \psi^i\wedge\omega_j
  -   \sum_\alpha\omega^{i}_{\alpha}\wedge\omega^\alpha_j   =  \frac{1}{2}\sum_{kl}
  R^{i}_{\ jkl}\omega^k\wedge\omega^l;
\eqno{(3.17)}$$
 $$d\omega^{\alpha}_{\beta}   +  \sum_k\omega^{\alpha}_{k}\wedge\omega^k_\beta
  =
  -   \sum_\alpha\omega^{\alpha}_{\alpha}\wedge\omega^\alpha_\beta
   =  \frac{1}{2}\sum_{kl}
  R^{\alpha}_{\ \beta
  kl}\omega^k\wedge\omega^l.
  \eqno{(3.18)}$$
   Denote
  $$ g_{_{ij}}=\langle Y_i,Y_j\rangle,  \quad
 g_{_{\beta\gamma}}=\langle \xi_\beta, \quad
 \xi_\gamma\rangle, \quad
 (g^{ij})=(g_{_{ij}})^{-1},\quad
 (g^{\beta\gamma})=(g_{_{\beta\gamma}})^{-1},$$
 $$
 R _{ij
  kl}= \sum _p g_{_{it}}
  R^{p }_{\ j
  kl},\quad
   R_{\alpha  \beta
  kl}= \sum_\nu g_{\alpha\nu}
  R^{\nu}_{\ \beta
  kl}.$$
 Then
the integrable conditions of the structure equations contain
  $$A_{ij, k}-A_{ik, j}=
  -  \sum_{\alpha\beta}g_{_{\alpha\beta}}
  (  B^\alpha_{ij}C^\beta_k  -  B^\alpha_{ik}C^\beta_j   ); \quad B^\alpha_{ij,
  k}-B^\alpha_{ik, j}=  g_{_{ij}}C^\alpha_k - g_{_{ik}}C^\alpha_j;
  \eqno{(3.19)}$$
  $$C^\alpha_{i, j}-C^\alpha_{j, i}
  =\sum_{kl}g^{kl}(B^\alpha_{ik}A_{lj}-B^\alpha_{jk}A_{li});
  \quad R_{\alpha\beta
  ij}  =  \sum_{kl\gamma\nu}g_{\alpha\gamma}g_{\beta\nu}  g^{kl}
(B^\gamma_{ik}B^\nu_{lj}-B^\nu_{ik}B^\gamma_{lj}); \eqno{(3.20 )}$$
  $$R_{ijkl}=\sum_{\alpha\beta}g_{_{\alpha\beta}}
(B^\alpha_{ik}B^\beta_{jl}-B^\alpha_{il}B^\beta_{jk})+
(g_{ik}A_{jl}-g_{il}A_{jk}) +(A_{ik}g_{jl}-
A_{il}g_{jk}).\eqno{(3.21)}$$
 Furthermore, we have
  $$\text{tr}(\mathbb A)=\frac{1}{2m}( m^2\kappa\pm1);\quad
  R_{ij}=\text{tr}(\mathbb
  A) g_{ij}+(m-2)A_{ij}-\sum_{kl\alpha\beta}
  g^{kl}g_{_{\alpha\beta}}B^\alpha_{ik}
  B^\beta_{lj} ;\eqno{(3.22)}$$
  $$(1-m)C^\alpha_i=\sum_{jk}g^{jk}B^\alpha_{ij,k};\quad
  \sum_{ijkl\alpha\beta}g^{ij}g^{kl}
 g_{_{\alpha\beta}}B^\alpha_{ik}B^\beta_{jl}=\frac{m-1}{m};\quad
 \sum_{ij } g^{ij}B^\alpha_{ij}=0,\forall\alpha. \eqno{(3.23)}$$
 From above we know that in case
$m\geq3$ all coefficients in the PDE system (3.10)-(3.12) are
determined by the conformal metric $g$,  the conformal second
fundamental form $\mathbb B$ and the normal connection
$\{\omega^\beta_\alpha\}$ in the conformal normal bundle $\mathbb
V$. Then we have

{\bf Theorem 3.3.} Two hypersurfaces $x:\mathbf M^m\rightarrow
\mathbb Q^{m+1}_p$ and $\widetilde{x}:\widetilde{\mathbf
M}^m\rightarrow \mathbb Q^{m+1}_p(m\geq3)$ are conformal equivalent
if and only if there exists a diffeomorphism $f:\mathbf
M\rightarrow\widetilde{ \mathbf M}$ which preserves the conformal
metric   and the conformal second fundamental form. In another word,
$\{g, \mathbb B\}$ is a complete invariants system of the
hypersurface $x:\mathbf M^m\rightarrow
\mathbb Q^{m+1}_p(m\geq3)$.

 Next we give the relations between the conformal
invariants induced above and $SO(n-p+1$,p+1)-invariants of $u:\mathbf
M\rightarrow\mathbb R^{n}_p$. We give also a conformal fundamental
theorem
for hypersurfaces  in $\mathbb R^n_p$.
  The psudo Euclidean space $\mathbb R^n_p$ has an inner product
  $\langle,\rangle$,  whose signature is $( \underbrace{+,\cdots,+}_{(n-p )-tiple},$
  $\underbrace{-,\cdots,-}_{  p -tiple})$. From the conformal
map
$$\sigma :\mathbb R^{n}_p\rightarrow{\mathbb Q}^n_p,\quad u\mapsto
[(\frac{  \langle u , u  \rangle  -1}{2},u, \frac{  \langle u , u  \rangle  +1}{2})],  \eqno{( 3.24)}$$
 we may recognize that
$\mathbb R^n_p\subset{\mathbb Q}^n_p$.
 Let  $u:\mathbf M\rightarrow \mathbb R^{n}_p$ be a submanifold.
 Let
 $\{e_1, \cdots , e_m\}$ be an local basis for $u$ with dual basis
 $\{\omega^1, \cdots , \omega^m\}$. Let $\{e_{m+1}, \cdots ,$ $ e_n\}$
 be a local basis of the normal bundle of $u$ in $\mathbb R^{n}_p$.
 Then we have the first and second fundamental forms $I, II$ and the
 mean curvature vector $\mathbf H$.  We may write
 $$I  =  \sum_{ij}I_{ij}\omega^i\otimes\omega^j,\quad
 II=\sum_{ij\alpha}h^\alpha_{ij}\omega^i\otimes\omega^je_\alpha$$
 $$(I^{ij}) = (I_{ij})^{-1},\quad
 \mathbf H = \frac{1}{m}\sum_{ij\alpha}I^{ij}
 h^\alpha_{ij}e_\alpha:=\sum_\alpha H^\alpha e_\alpha.$$
  Denote
 $\Delta_{\bf M}$ the Laplacian and $\kappa_{\bf M}$ the normalized
 scalar curvature for $I$.  It is easy to see that
 $$\Delta_{\bf M} u=m\mathbf H,
 \quad\kappa_{\bf M}=\frac{1}{m(m-1)}(m^2|\mathbf H|^2
 -|II|^2), $$
 where
   $$|\mathbf H|^2=\sum_{\alpha\beta}I_{\alpha\beta}H^\alpha H^\beta,
   I_{\alpha\beta}=(e_\alpha,e_\beta);\quad
    |II|^2=\sum_{ijkl\alpha\beta}
   I_{\alpha\beta}  I^{ik}  I^{jl}h^\alpha_{ij}h^{\beta}_{kl}.$$
 In fact, from the structure equations
 $$\mathrm{d}u=\sum_i\omega^iu_i,
  \quad\mathrm{d}u_i=\sum_j\theta^j_iu_j+\sum_\alpha\theta^\alpha_ie_\alpha,
  \quad\mathrm{d} e_\alpha
  =\sum_j\theta^j_\alpha u_j+\sum_\beta\theta^\beta_\alpha e_\beta,
  \eqno{(3.26)}$$
  we have
  $$\sum_ju_{i,j}\omega^j=\mathrm{d}u_i-\sum_j\theta^j_iu_j
  =\sum_\alpha\theta^\alpha_ie_\alpha,
 \quad u_{i,j} =\sum_\alpha h^\alpha_{ij}e_\alpha.\eqno{(3.27)}$$
 For $x=\sigma\circ u:M\rightarrow \mathbb{R}^n_p$, there is a global lift
 $$y:\mathbf M\rightarrow C^{n+1},\quad y=( \frac{  \langle u , u  \rangle  -1}{2},u, \frac{  \langle u , u  \rangle  +1}{2})
.$$ So we will get
  $$\langle\text dy, \text dy\rangle
  =  \langle \text du , \text du  \rangle  =I;
  \quad \Delta=\Delta_{\bf M};\quad \kappa=\kappa_{\bf M}.
 \eqno{( 3.28)}$$
  It follows from (3.25) that
  $$\langle\Delta Y, \Delta Y\rangle-m^2\kappa
  =\frac{m}{m-1}(|II|^2-m|\mathbf H|^2).
  \eqno{(3.29 )}$$
  Therefore the conformal metric of $x$
  $$g  = \pm \frac{m}{m-1}(|II|^2-m|\mathbf H|^2)
    \langle \text du , \text du  \rangle  :=e^{2\tau}I.
  \eqno{(3.30)}$$
 Let
 $$y_i = e_i (y ) = ( 0 ,u_i,   0 ) +  (u, u_i)  ( 1,\mathbf 0, 1),
 \zeta_\alpha = ( 0 ,e_\alpha,  0 )  + (u, e_\alpha) ( 1,\mathbf 0, 1).$$
   Through some
  direct calculation it reaches
  $$Y=e^\tau y,\quad Y_i = e_i ( Y) =e^\tau(\tau_iy+ y_i),\quad
  \xi_\alpha=H_\alpha y+\zeta_\alpha,\eqno(3.31)$$
  $$-e^\tau
  N=\frac{1}{2}(|\nabla \tau|^2+|\mathbf H|^2)
 y+\sum_i\tau^i
  y_i+\sum_\alpha H^\alpha\zeta_\alpha+( 1,\mathbf 0, 1),\eqno(3.32)$$
  where $\tau^i=\sum_{j}I^{ij}\tau_j,(I^{ij})=(I_{ij})^{-1};\quad
  |\nabla\tau|^2=\sum_{i}\tau_i\tau^i;\quad
  H_\alpha=\sum_{\beta}I_{\alpha\beta}
   H^\beta.$

 By a direct calculation we get the following expression of
the conformal invariants $\mathbb A,\mathbb B,$ and $\Phi$:
  $$    A_{ij} = \tau_i\tau_j -
  \sum_{\alpha}h^\alpha_{ij}H_\alpha -\tau_{i,
  j}  -  \frac{1}{2}(|\nabla\tau|^2+|\mathbf H|^2) I_{ij}, \eqno{(3.33)}$$
  $$  B^\alpha_{ij}  = e^{ \tau} ( h^\alpha_{ij}-H^\alpha I_{ij} ) ,
  \quad e^{ \tau} C^\alpha_i=  H^\alpha\tau_i
 -\sum_{j}h^\alpha_{ij}\tau^j-H^\alpha_{, i}
 ,\eqno{(3.34)}$$
 where $\tau_{i, j}$ is the
Hessian of $\tau$ respect to $I$ and $H^\alpha_{, i}$ is the covariant
derivative of the mean curvature vector field of $u$ in the normal
bundle $N(\mathbf M)$ respect to $I$.

  Now we consider the case that $u:\mathbf M\rightarrow\mathbb
R^n_p$ is a hypersurface.
Observing the PDE system (3.10)-(3.12), from   Theorem 3.3  we have

 {\bf Theorem 3.4. } Two hypersurfaces $u,\tilde u:\mathbf
 M\rightarrow\mathbb R^n_p(n\geq4)$ are conformally equivalent if and
 only if there exists a diffeomorphism $f:\mathbf M\rightarrow\mathbf
 M$ which preserves the conformal metric   and the conformal second fundamental form
 $\{g, \mathbb B\}$.

  {\bf Remark 3.1.} For psudo sphere   space with index $p$
  $$\mathbb
  S^n_p =  \{  u = ( u^1,\cdots,u^{n+1} )\in \mathbb{R}^{n+1}|
  \ \ \ \ \ \ \ \ \ \ \ \ \ \ \ \
  \ \ \ \ \ \  \ \ \ \ \ \ \ \ $$
  $$ \ \ \ \ \ \ \ \ \ \ \ \ \ \ \ \
  \ \ \ \ \ \  \langle  u ,  u  \rangle  :=
    (u^1)^2  +  \cdots +  ( u^{n-p +1 })^2 - ( u^{n-p+2} )^2  -  \cdots -  ( u^{n+1} )^2=1\}
  $$
   and   psudo hyperbolic  space with index $p$
  $$\mathbb H^n_p=  \{  u = ( u^1,\cdots,u^{n+1} )
  \in \mathbb{R}^{n+1}|\ \ \ \ \ \ \ \ \ \ \ \ \ \ \ \
  \ \ \ \ \ \ \ \ \ \ \ \ \ \ $$
  $$ \ \ \ \ \ \ \ \ \ \ \ \ \ \ \ \
  \ \ \ \ \ \
  \langle  u ,  u  \rangle  :=     (u^1)^2  +  \cdots +  ( u^{n-p   })^2
    - ( u^{n-p+  1} )^2  -  \cdots -  ( u^{n+1} )^2  =-1\}$$
     we obtain analogous conclusion:
      $$    A_{ij} = \tau_i\tau_j -
  \sum_{\alpha}h^\alpha_{ij}H_\alpha -\tau_{i,
  j}  -  \frac{1}{2}(|\nabla\tau|^2+|\mathbf H|^2 -\epsilon ) I_{ij}, \eqno{(3.35)}$$
  $$  B^\alpha_{ij}  = e^{ \tau} ( h^\alpha_{ij}-H^\alpha I_{ij} ) ,
  \quad e^{ \tau} C^\alpha_i=  H^\alpha\tau_i
 -\sum_{j}h^\alpha_{ij}\tau^j-H^\alpha_{, i}
 ,\eqno{(3.36)}$$
 where $\epsilon$ corresponds the sectional curvature of
 psudo sphere   space or  psudo hyperbolic  space with index $p$.\\

\par\noindent
{\bf {\S} 4. The first variation of the conformal volume functional}
\par\medskip
  Let $x_0:\mathbf M\rightarrow {\mathbb Q}^n_p$ be a compact oriented regular
  submanifold with boundary $\partial \mathbf M$. Suppose that local
   basis $\{e_1,\cdots,e_m\}$ on $\mathbf M$ satisfy the orientation. Denote $g_{ij}=g(e_i,e_j)$.
  If the conformal metric $g$ has $s$ negative signature and  $( g_{ij} )
  =  ( - I_s) \oplus  ( I_{m-s} )$, we call $\{e_1,\cdots,e_m\}$
   a local orthonormal basis for g.
   In the following let $\{e_1,\cdots,e_m\}$ be
   a local  orthonormal basis for $g$ with dual basis $\{\omega^1,\cdots,\omega^m\}$.

   We define the
  generalized Willmore functional $\mathbb{W}(\mathbf M)$ as the
  volume functional of the conformal metric $g$:
  $$\mathbb{W}(\mathbf M)={\text {Vol}}_g(\mathbf M)=\int_{\bf M}\text d {\bf M}_g.$$
   The conformal volume element $\text d {\bf M}_g$ is defined by
   $$\text d {\bf M}_g  =
   \omega^1\wedge\cdots \wedge \omega^m,$$
   which is well-defined.

   Let
  $x:\mathbf M\times  \mathbb R\rightarrow {\mathbb Q}^n_p$ be an
  admissible variation of $x_0$ such that $ x(\cdot,t)=x_t$ and $
  dx_t(\mathrm{T}_p\mathbf M)= dx_0(\mathrm{T}_p\mathbf M)$ on
  $\partial \mathbf M$ for each small $t$.
 For each $t$, $x_t$ has  the conformal metric $g_{_t}$.
  As  in \S 3, we have a moving frame $\{Y,N,Y_i,\xi_\alpha\}$ in
  $\mathbb R^{n+2}_{p+1}$ along $\mathbf M\times \mathbb R$ and the
  conformal volume $W(t)=\mathbb{W}(x_t)$. Let $\{\xi_\alpha\}$ be a local
   orthonormal  basis for the conformal normal bundle $\mathbb V_t$ of $x_t$. Denote
  $\widetilde{\mathrm{d}}$ and $\mathrm{d} $ the differential operators
  on $\mathbf M\times \mathbb R$ and $\mathbf M$, respectively. Then
  we have
  $$\widetilde{\mathrm{d}} =  \mathrm{d}
  +  \text dt\wedge\frac{\partial}{\partial t}\eqno{(4.1)}$$
  on $\mathrm{T}^*(\mathbf M\times \mathbb R)=\mathrm{T}^*\mathbf
  M\oplus \mathrm{T}^*\mathbb R$.
 We also have
 $$ \mathrm{d} \circ\frac{\partial}{\partial t}
 =  \frac{\partial}{\partial t} \circ              \mathrm{d} .
 \eqno{(4.2)}$$
  Denote
  $ P =  ( Y,N,Y_i,\xi_\alpha )^T$. Suppose that
  $  \mathrm{d} P= \Omega P,
  \frac{\partial}{\partial t}P = L P,$
 where
 $$ \Omega= \left(
         \begin{array}{cccc}
  0   &  0    &    \omega^j  &  \mathbf{0} \\
  0   &  0    &  \psi^j  &   \phi^\beta\\
  -\psi_i    & -\omega_i    &   \omega_i^j  &   \omega_i^\beta\\
   -\phi_\alpha    & \mathbf{0} &   \omega_\alpha^j  &   \omega_\alpha^\beta
         \end{array}
 \right),  L =
 \left(
         \begin{array}{cccc}
  w   &  0    &    v^j  &    v^\beta \\
  0   &  -w    &   u^j  &    u^\beta\\
  -u_i    & -v_i    &   L_i^j  &     L_i^\beta\\
   -u_\alpha    & -v_\alpha  &     L_\alpha^j  &    L_\alpha^\beta
         \end{array}
 \right)
 .
  $$
   From (4.2) it is easy to get
  $$ \frac{\partial}{\partial t}
  \Omega = \mathrm{d}  L  +  L \Omega  -
  \Omega L. \eqno{(4.3)}$$
  Therefore we have
  $$\frac{\partial\omega^i}{\partial
  t}=\sum_j    (       v^i_{,j}      +      L_{  j  }^{  i  }
  -       \sum_{k\alpha\beta}  g_{\alpha\beta}v^\alpha
  B^\beta_{kj}   g^{ik}
  )\omega^j   +  \sum_\alpha  v^\alpha  \omega^i_\alpha
  +  w\omega^i,\quad
  L^\alpha_i=v^\alpha_{,i}+\sum_jB^\alpha_{ij}v^j,\eqno{(4.4)}$$
 where $\{  v^i_{,j}  \}$ is
 the covariant derivative of $\sum v^i e_i$ with respect to
 $g_{_t}$ and $\{v^\alpha_{,i}\}$ is
 the  covariant
 derivative of $\sum v^\alpha \xi_\alpha$.
 Here we have used the notations of conformal invariants
 $\{A_{ij},B^\alpha_{ij},C^\alpha_i\}$ for  $x_t$ defined in \S 3.
 Furthermore we have
  $$\frac{\partial\omega^\alpha_i}{\partial
  t}=\sum_j(L^\alpha_{i,j}+\sum_kL_{i}^{k}B^\alpha_{kj} -\sum_\beta
  B^\beta_{ij}L^\alpha_\beta
  +A_{ij}v^\alpha-v_i C^\alpha_j)\omega^j+u^\alpha
  \omega_i,\eqno{(4.5)}$$
  where   $\{L^\alpha_{i,j}\}$ is
 the  covariant
 derivative of
 $\sum_{i\alpha}L^\alpha_i\omega^i\xi_\alpha$.
  Using
 (4.4) and (4.5)   we get
 $$\frac{\partial B^\alpha_{ij}}{\partial t}
 +   w   B^\alpha_{ij}       =
 v^\alpha_{,ij}   +A_{ij}v^\alpha     +
 \sum_{kl\gamma}  g^{kl}  B^\alpha_{ik}B^\gamma_{lj}v_\gamma
 \ \ \ \ \ \ \ \ \ \ \ \ \ \ \ \ \ \ \ \ \ \ \ \ \ \  \   \ \ \ \ \ \ \ \ \ \ \
 \ \ \ \ \ \ \ $$
 $$ \ \ \ \ \ \ \ \ \ \ \ \ \ \ \ \ \ \ \ \ \ \ +   u^\alpha g_{ij} +      \sum_k  L_{i}^{k}B^\alpha_{kj}  -  \sum_\gamma
 B^\gamma_{ij}L^\alpha_\gamma
  +   \sum_k  v^{k}B^\alpha_{ik,j} -        v_i C^\alpha_j.
 \eqno{(4.6)}$$
 It follows from (3.19) and (3.23) that
 $$\frac{m-1}{m}w    =
 \sum_{ijkl\alpha\beta}
  g_{_{ \alpha \beta }} g^{ik} g^{jl} B^\beta_{kl} (
  v^\alpha_{,ij}        +        A_{ij} v^\alpha     +
 \sum_{kl\gamma}  g^{kl}  B^\alpha_{ik}B^\gamma_{lj}v_\gamma ).
 \eqno{(4.7)}$$
  Now we calculate the first variation
 of the conformal volume functional
 $$W(t)={\text {vol}}(g_{_t})
 =\int_{\bf M} \omega^1\wedge\cdots\wedge\omega^m=\int_{\bf M} \text
 d {\bf M}_g,$$
 where $\text d {\bf M}_g $ is the volume
 for $g_{_t}$. From (4.4)   we get
 $$W'(t)=\sum_i\int_{\bf M} \omega^1\wedge\cdots\wedge\frac{\partial\omega^i}{\partial
 t}\wedge\cdots\wedge\omega^m
 =   \int_{\bf M} \sum_i v^i_{,i} \text
 d {\bf M}_g  +    m\int_{\bf M}w \text d {\bf M}_g.\eqno{(4.8)}$$
 From the fact that
the variation is admissible we know $v^i =0,v^\alpha=0$ and
$v^\alpha_{,i}=0$ on $\partial \mathbf M$. It follows from (4.8)
and Green's formula that
$$W'(t)=\frac{m^2}{m-1}\int_{\bf M}\sum_\alpha
 v^\alpha \quad[  \quad
 \sum_{ijkl\beta}
  g_{_{ \alpha \beta }} g^{ik} g^{jl}\cdot \ \ \ \ \ \ \ \ \ \ \ \ \ \ \ \ \ \ \ \
  \ \ \ \ \ \ \ \ \ \ \ \ \ \ \ \ \ \ \ \ \ \ \ \ $$
  $$\ \ \ \ \ \ \ \ \ \ \ \ \ \ \ \ \ \ \ \ \ \ \ \
  \ \ \ \ \ \ \ \ \ \ \ \ \ \ \ \ \ \ \ \ \ \ \ \ \cdot(
  B^\beta _{ij,kl}         +        A_{ij}B^\beta_{kl}        +
 \sum_{ r q \gamma  \nu  }  g_{_{ \gamma  \nu  }}  g^{rq}
 B^\beta_{ir}           B^\gamma_{qj}           B^\nu_{kl}   )
 \quad]   \text d {\bf M}_g.\eqno{(4.9)}$$

 It follows from  (4.9) that

 {\bf Theorem 4.1.} The variation of the conformal volume functional
 depends only on the normal component of the variation field
 $\frac{\partial Y}{\partial t}$.  A submanifold $x:
 \mathbf{M}\rightarrow \mathbb{Q}^n_p$
is a Willmore  submanifold ({\sl i.e.}, a critical submanifold to the
conformal volume functional) if and only if
  $$\sum_{ijkl\beta}
  g_{_{ \alpha \beta }} g^{ik} g^{jl} (
  B^\beta _{ij,kl}         +        A_{ij}B^\beta_{kl}          +
 \sum_{ r q \gamma  \nu  }  g_{_{ \gamma  \nu  }}  g^{rq}
 B^\beta_{ir}           B^\gamma_{qj}           B^\nu_{kl}    )
  =0,\quad
\forall\alpha .\eqno{(4.10)}$$

  We call the equation (4.10) the
 Euler-Lagrange equations   or  Willmore equations.
 Using (3.22) and (3.23) we can write the  Willmore equations  (4.10) as
  $$\sum_{ \beta}
  g_{_{ \alpha \beta }}   \big[   \sum_{ij} g^{ij} C^\beta_{i,j}  +
  \sum_{ijkl} g^{ik} g^{jl} (\frac{1}{m-1}R_{ij}-A_{ij})B^\beta_{ kl }  \big]
  =0
,\quad \forall\alpha .\eqno{(4.11)}$$

 {\bf Theorem 4.2.}  Any stationary ({\sl means that whose curvature vector
 is vanishing})
regular surface in psudo Euclidean space $\mathbb R^n_p$,
psudo sphere space $\mathbb S^n_p$ and
psudo hyperbolic  space $\mathbb H^n_p$ is Willmore.

 {\sl Proof\ \ } Let
$u:\mathbf{M}\rightarrow \mathbb R^n_p$ be a regular surface, whether space-like or time-like. Let
$\{e_1,e_2\}$ be a local   basis of $\langle\text du,\text du\rangle$
and $\{ e_\alpha\}_{\alpha=3}^n$ a local basis for the normal bundle.
  If $x$ is a stationary regular
surface, we have $H^\alpha\equiv0,\forall\alpha$.
  From (3.33) and (3.34) we get
  $$\sum_{ijkl}g^{ik} g^{jl} A_{ij}B^\beta_{kl}=
  \sum_{ijkl}g^{ik} g^{jl}  B^\beta_{kl}(\tau_i\tau_j
   -\tau_{i,j})=   e^{-3\tau}
  \sum_{ijkl}  I^{ik} I^{jl}  h^\beta_{kl}(\tau_i\tau_j
   -\tau_{i,j})
       .\eqno{(4.12)}$$
       Now we know from (3.34)  that
       $$   -  e^\tau C^\beta_{i }  =
  \sum_{  kl}  I^{ k l}  h^\beta_{ ik } \tau_l
  := W^\beta_i
  .\eqno{(4.13)}$$
  From (3.14)  we have
  $$\sum_j e^\tau   C^\beta_{i,j}  \omega^j   =
  \mathrm d (  e^\tau   C^\beta_{i }  )  -
      e^\tau   C^\beta_{i }  \mathrm d\tau
      +  \sum_\gamma   e^\tau   C^\gamma_{i }  \theta_\gamma^\beta
      -  \sum_k   e^\tau   C^\beta_k  \omega_i^k\ \ \ \ \ \ \ \
      \ \ \ \ $$
      $$  =  - \mathrm d   W^\beta_{i }
      +  W^\beta_{i } \mathrm d \tau
      -  \sum_\gamma W^\gamma_{i }  \theta_\gamma^\beta
      +  \sum_k   W^\beta_k  \omega_i^k
      .\eqno{(4.14)}$$
        Combining with
  $$ \omega^k_i  =  \theta^k_i  +  \tau^k \sum_j  I_{ij}\omega^j
  -  \tau_i\omega^k  +  \delta^k_i \mathrm d\tau$$
  and (4.14) we get
  $$   e^\tau   C^\beta_{i,j}  =
  2 W^\beta_{i } \tau_j  +
   W^\beta_{j }  \tau_i  -
  \sum_k  W^\beta_{ k }  \tau^k  I_{ij}  -
     W^\beta_{i,j }
       ,\eqno{(4.15)}$$
       where $ W^\beta_{i,j }$ is the covariant differential of $ W^\beta_{i  }$ with
       respect to the first fundamental form $I$ of $u$.
       Therefore
       $$ \sum_{ij} g^{ij} C^\beta_{i,j}
       =   e^{-3\tau}
  \sum_{ijkl}  I^{ik} I^{jl}  h^\beta_{kl}(\tau_i\tau_j
   -\tau_{i,j})
       .\eqno{(4.16)}$$
  Whether the regular surface $u$ is space-like or time-like,
  if we choose $\{e_1,e_2\}$ orthonormal, then
  a direct calculation leads to
  $$ \sum_{ijkl}g^{ik} g^{jl} R_{ij}B^\beta_{kl}
  =  0.\eqno{(4.17)}$$
   Thus we have (4.11) from (4.12), (4.16) and (4.17), which implies that $u$ is Willmore.

 One can verify
  that stationary regular surfaces in
  $\mathbb S^n_p$ and $\mathbb H^n_p$ are also Willmore. $\Box$

  {\bf Remark 4.1.} In some conferences, a surface in psudo Riemannian space
  forms with vanishing mean curvature vector is also called {\sl
  maximal} or {\sl minimal}. But in this time the volume functional of
  the surface is not really {\sl maximal} or {\sl minimal}. So we take the
  place of the above two terms by {\sl stataionary}
  ({\sl also see} [1]).\\

\par\noindent
{\bf {\S} 5. Conformal isotropic submanifolds in ${\mathbb Q}^n_p$}
\par\medskip

{\bf Definition 5.1.} We call an m-dimensional submanifold $x:\mathbf
M\rightarrow {\mathbb Q}^n_p$ is conformal isotropic if there exists
a smooth function $\lambda$ on $\mathbf M$ such that
$$\mathbb A+\lambda g\equiv0\ \text{ and}\ \ \Phi\equiv0.\eqno{(5.1)}$$
\indent
 From precious discuss in  \S 3
 we can easily verify

{\bf Proposition 5.1.} If $u:\mathbf M\rightarrow \mathbb R^n_p$ is a
stationary regular submanifold with constant scalar curvature,
then $x=\sigma \circ u$ is a conformal isotropic submanifold in
${\mathbb Q}^n_p$.

{\bf Remark 5.1.} The same conclusion holds on $\mathbb S^n_p$ or
$\mathbb H^n_p$.

Suppose that $x:\mathbf M\rightarrow {\mathbb Q}^n_p$ is a conformal
isotropic submanifold. Then we get
$$\mathrm{d}N+\lambda \mathrm{d}Y=0, \quad \mathrm{d}\lambda\wedge \mathrm{d}Y
=\sum_{i=1}^m(\mathrm{d}\lambda\wedge\omega^i)Y_i=0. \eqno{(5.2)}$$

 Since
$\{Y_1, \cdots , Y_m\}$ are linearly independent,
$$ \mathrm{d}\lambda\wedge \omega^i=\sum_{j=1}^mE_j(\lambda)\omega^j\wedge\omega^i =0.
\eqno{(5.3)}$$

If $\mathbf M$ is connected, we get
$$\lambda=\text{constant}, \eqno{(5.4)}$$
which concludes from (2.3) that
$$\kappa=\text{constant}.\eqno{(5.5)}$$

By (5.2) we can find a constant vector $\mathbf c\in \mathbb R^{n+2}_{p+1}$ such that
$$N+\lambda Y=\mathbf c. \eqno{(5.6)}$$

It follows that
$$\langle Y, \mathbf c\rangle=1, \langle\mathbf c, \mathbf c\rangle=2\lambda=\text{constant}. \eqno{(5.7)}$$

Then we look into three cases.

Case 1: $\langle\mathbf c, \mathbf c\rangle=0$.  By making use of a
Lorenzian rotation in $\mathbb R^{n+2}_{p+1}$ when necessary,  we may
assume that
$$\mathbf c=(-1,0, -1). \eqno{(5.8)}$$
Letting
$$Y=(x_p,u,x_{n+2}), \eqno{(5.9)}$$
by (5.7) and $Y\in C^{n+1}$ we have
$$Y=(\frac{  \langle u , u  \rangle  -1}{2},u, \frac{  \langle u , u  \rangle  +1}{2}). \eqno{(5.10)}$$
Then $x$ determines a submanifold $u:\mathbf M\rightarrow\mathbb
R^{n}_p$ with $I=  \langle \text du , \text du  \rangle  =\langle\text dy, \text
dy\rangle=g$, which implies that
$$\kappa_{\bf M}=\kappa=\text{constant}. \eqno{(5.11)}$$
From (5.7) and (2.13 ) we have $H^\alpha=0$,  {\sl i.e.},  $u$ is a
stationary submanifold in $\mathbb R^n_p$.  In this case $x$ is
conformal equivalent to the image of $\sigma $ of a stationary
submanifold with
constant scalar curvature in $\mathbb R^n_p$.

Case 2: $\langle\mathbf c, \mathbf c\rangle=- r^2, r>0. $ By making use
of a Lorenzian rotation in $\mathbb R^{n+2}_{p+1}$ when necessary,  we
may assume that
$$\mathbf c=(\mathbf{0},r). \eqno{(5.12)}$$
 Letting
$$Y=(u/r,x_{n+2}), \eqno{(5.13)}$$
by (5.7)  we have
$$x_{n+2}=1/r. \eqno{(5.14)}$$
So
$$Y=(u,1)/r,   \langle u , u  \rangle  = 1 . \eqno{(5.15)}$$
Then $x$ determines a submanifold $u:\mathbf M\rightarrow \mathbb
S^n_p$ with $I/r^2=  \langle \text du , \text du  \rangle  /r^2=\langle\text dy, \text
dy\rangle=g$, which implies that $\kappa_{\bf
M}=\kappa/r^2=\text{constant}. $ From (5.7) and (2.21) we have
$H^\alpha=0$,  {\sl i.e.}, $u$ is a stationary submanifold in $\mathbb
S^n_p$. In this case $x$ is conformal equivalent to the image of
$\sigma_+$ of a stationary submanifold with constant scalar
curvature in
$\mathbb S^n_p$.

Case 3: $\langle\mathbf c, \mathbf c\rangle=r^2,  r>0. $ By making use
of a Lorenzian rotation in $\mathbb R^{n+2}_{p+1}$ when necessary, we
may assume that
$$\mathbf c=(- r,\mathbf{0}). \eqno{(5.16)}$$
Letting
$$Y=(x_p,u/r), \eqno{(5.17)}$$
by (5.7)  we have
$$x_p=1/r. \eqno{(5.18)}$$
So
$$Y=(1,u)/r,   \langle u , u  \rangle  =-1 . \eqno{(5.19)}$$
Then $x$ determines a submanifold $u:\mathbf M\rightarrow \mathbb
H^n_p$ with $I/r^2=  \langle \text du , \text du  \rangle  /r^2=\langle\text dy, \text
dy\rangle=g$, which implies that $\kappa_{\bf
M}=\kappa/r^2=\text{constant}. $ From (5.7) and (2.26) we have
$H^\alpha=0$, {\sl i.e.}, $u$ is a stationary submanifold in $\mathbb
H^n_p$. In this case $x$ is conformal equivalent to the image of
$\sigma_{- }$ of a stationary submanifold with constant scalar
curvature in
$\mathbb H^n_p$.

 So combining with   Proposition 5.1
and   Remark 5.1 we get

 {\bf
Theorem 5.2.} Any conformal isotropic submanifold in ${\mathbb Q}^n_p$ is conformal equivalent to a stationary submanifold with
constant scalar
curvature in $\mathbb R^n_p,\mathbb S^n_p$, or $\mathbb H^n_p$.\\

{\bf Acknowledgements: }
 The authors would like to express his
gratitude to Professor Changping   Wang  for his warm-hearted inspiration and
the support of BICMR.

\vskip 1cm
\begin{flushleft}

Changxiong Nie\\
\textsl{ Faculty of Mathematics and Computer Sciences  }\\
\textsl{ Hubei University  }\\
\textsl{ 430062 Wuhan  }
 \\
\textsl{ People's Republic of China  }\\
\texttt{  chxnie@163.com}

\end{flushleft}

\end{document}